\documentclass[12pt]{article}
\usepackage{mathrsfs}
\usepackage{amssymb}
\usepackage{amsmath}
\usepackage{graphicx}
\usepackage{amsmath,amsfonts}

\newtheorem{theorem}{Theorem}[section]

\newcommand{\be}{\begin{equation}} \newcommand{\ee}{\end{equation}}
\newcommand{\ba}{\begin{align}} \newcommand{\ea}{\end{align}}
\newcommand{\baa}{\begin{align*}} \newcommand{\eaa}{\end{align*}}
\newcommand{\ben}{\begin{enumerate}} \newcommand{\een}{\end{enumerate}}
\newcommand{\bi}{\begin{itemize}} \newcommand{\ei}{\end{itemize}}

\newcounter{example}

\begin{document}
\title{Decomposition of supercritical linear-fractional branching processes} % insert title - use \\ if it requires more than one line.

\author{Serik Sagitov\footnote{Mathematical Sciences, Chalmers University of Technology and University of Gothenburg}
 \ and Altynay Shaimerdenova\footnote{Al-Farabi Kazakh National University}} 
 \maketitle

\begin{abstract}
It is well known that a supercritical single-type Bienyam\'e-Galton-Watson process can be viewed as a decomposable branching process formed by two subtypes of particles: those having infinite line of descent and those who have finite number of descendants. 
In this paper we analyze such a decomposition for the  linear-fractional Bienyam\'e-Galton-Watson  processes with countably many types.
\end{abstract}


\section{Introduction}

The Bienaym\'e-Galton-Watson (BGW-) process is a basic model for the stochastic dynamics of the size of a population formed by independently reproducing particles. It has a long history \cite{HS} with its origin dating back to 1837. This paper is devoted to the BGW-processes with countably many types.
One of the founders of the theory of multi-type branching processes is B.A. Sevastyanov  \cite{Sev}, \cite{Sew}.

A single-type BGW-process is a Markov chain $\{Z^{(n)}\}_{n=0}^\infty$ with countably many states $\{0,1,2,\ldots\}$. The evolution of the process is described by a probability generating function
\be
f(s)=\sum_{k=0}^\infty p_ks^k, \qquad p_1<1,
\label{fs}
\ee
where $p_k$ stands for the probability that a single particle produces exactly $k$ offspring. If particles reproduce independently with the same reproduction law \eqref{fs}, then the chain $\{Z^{(n)}\}_{n=0}^\infty$ represents consecutive generation sizes. In this paper, if not specified otherwise, we assume that  $Z^{(0)}=1$, the branching process stems from a single particle.  Due to the reproductive  independence it follows that
$f^{(n)}(s)=\mathbb E(s^{Z^{(n)}})$ is the $n$-th iteration of $f(s)$.

Since zero is an absorbing state of the BGW-process,
$q^{(n)}=\mathbb P(Z^{(n)}=0)$ monotonely increases to a limit $q$ called the extinction probability. The latter is implicitly determined as a minimal non-negative solution of the equation
\be
f(x)=x.
\label{fx}
\ee
A key characteristic of the BGW-process is the mean offspring number $M=f'(1)$. In the subcritical ($M<1$) and critical ($M=1$) cases the process is bound to go extinct $q=1$,
while in the supercritical case ($M>1$) we have $q<1$. Clearly $q=0$ if and only if $p_0=0$.

In the supercritical case the number of descendants of the progenitor particle is either finite with probability $q$ or infinite with probability $1-q$. Recognizing that the same is true for any particle appearing in the BGW-process we can distinguish between  {\it skeleton particles} having an infinite line of descent  \cite{O} and {\it doomed particles} having a finite line of descent. Graphically we get a picture of the genealogical tree similar to that given in Figure 1. 

If we disregard the doomed particles, the skeleton particles form a BGW-process with a transformed reproduction law excluding extinction
\be
\tilde f(s)={f(s(1-q)+q)-q\over 1-q}\label{HST}
\ee
and having the same mean $\tilde M=M>1$. Formula \eqref{HST} is usually called the Harris-Sevastyanov transformation.  On the other hand, the doomed particles form another branching process corresponding to the supercritical branching process conditioned on extinction. The doomed particles produce only doomed particles according to another transformation of the reproduction law $\hat f(s)={f(sq)/q}$, which is usually called the dual reproduction law and has mean $\hat M=f'(q)<1$.
The supercritical BGW-process as a whole can be viewed as a decomposable branching process with two subtypes of particles \cite[Ch. 1.12]{AN}.
Each skeleton particle must produce at least one new skeleton particle and also can give rise to a number of doomed particles. 
In Section \ref{sHS} we describe in detail this decomposition for the single type supercritical BGW-processes.

\begin{figure}
\centering
\includegraphics[height=3.5cm]{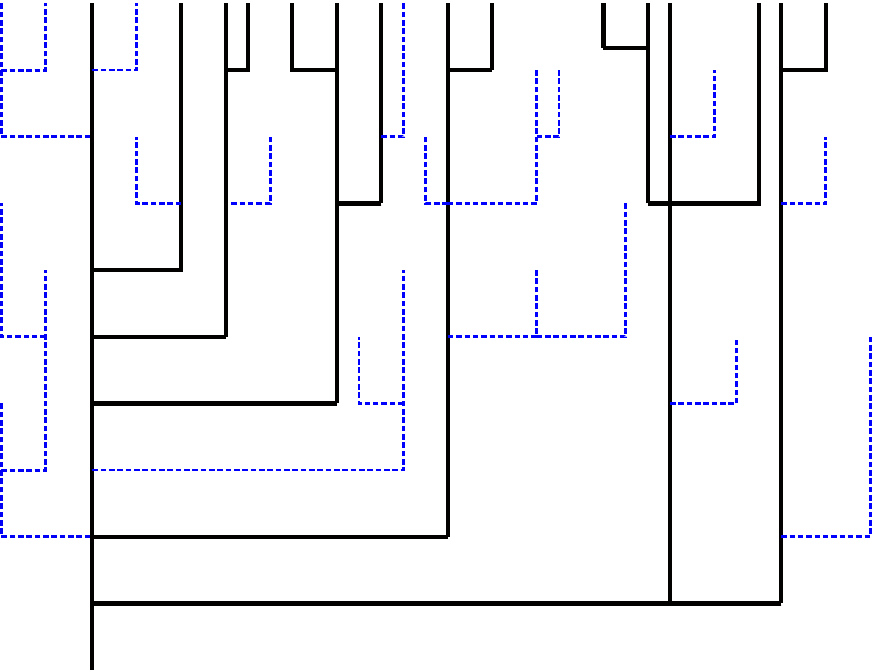}
\caption{An example of a BGW-tree up to level $n=10$. Solid lines represent the infinite lines of descent and dotted lines represent the finite lines of descent.}
\label{fig1}
\end{figure}

In the special case when the reproduction generating function \eqref{fs} is linear-fractional many characteristics of the BGW-process can be computed in an explicit form  \cite{KRS}. In Section \ref{sLF} we summarize explicit results concerning decomposition of a supercritical single-type BGW-processes.

Section \ref{sC} presents the BGW-processes with countably many types. Our focus is on the linear-fractional case recently studied in \cite{S}. 
The main results of this paper are collected in Section \ref{sMR} and their derivation is given in Section \ref{sP}. The remarkable fact that a supercritical branching process conditioned on extinction is again a branching process was recently  established in \cite{JL}  in a very general setting. In general, the transformed reproduction laws are characterized in an implicit way and are difficult to analyse. This paper presents a case where the properties of the skeleton and doomed particles are very transparent.

\section{Decomposition of a supercritical single-type BGW-process}\label{sHS}

The BGW-process is a time homogeneous Markov chain  with transition probabilities satisfying
\begin{align*}
\sum_{j=0}^\infty P^{(n)}_{ij}s^j=\left(f^{(n)}(s)\right)^i.
\end{align*}
In the supercritical case with mean $M>1$ and extinction probability $q<1$, using the property
$\sum_{j=0}^\infty P^{(n)}_{ij}q^j=q^i$,
we can get another set of transition probabilities putting
\begin{align*}
\hat P^{(n)}_{ij}:=P^{(n)}_{ij}q^{j-i}.
\end{align*}
The transformed  transition probabilities also possess the branching property
\begin{align*}
\sum_{j=0}^\infty \hat P^{(n)}_{ij}s^j=\left(\hat f^{(n)}(s)\right)^i,
\end{align*}
where $\hat f^{(n)}(s)={f^{(n)}(sq)/ q}$ is the $n$-th iteration of the so-called dual generating function
\[\hat f(s)={f(sq)\over q}=\sum_{k=0}^\infty \hat p_ks^k, \qquad \hat p_k=p_kq^{k-1}, \ k\ge0.\]

\begin{figure}
\centering
\includegraphics[width=10cm]{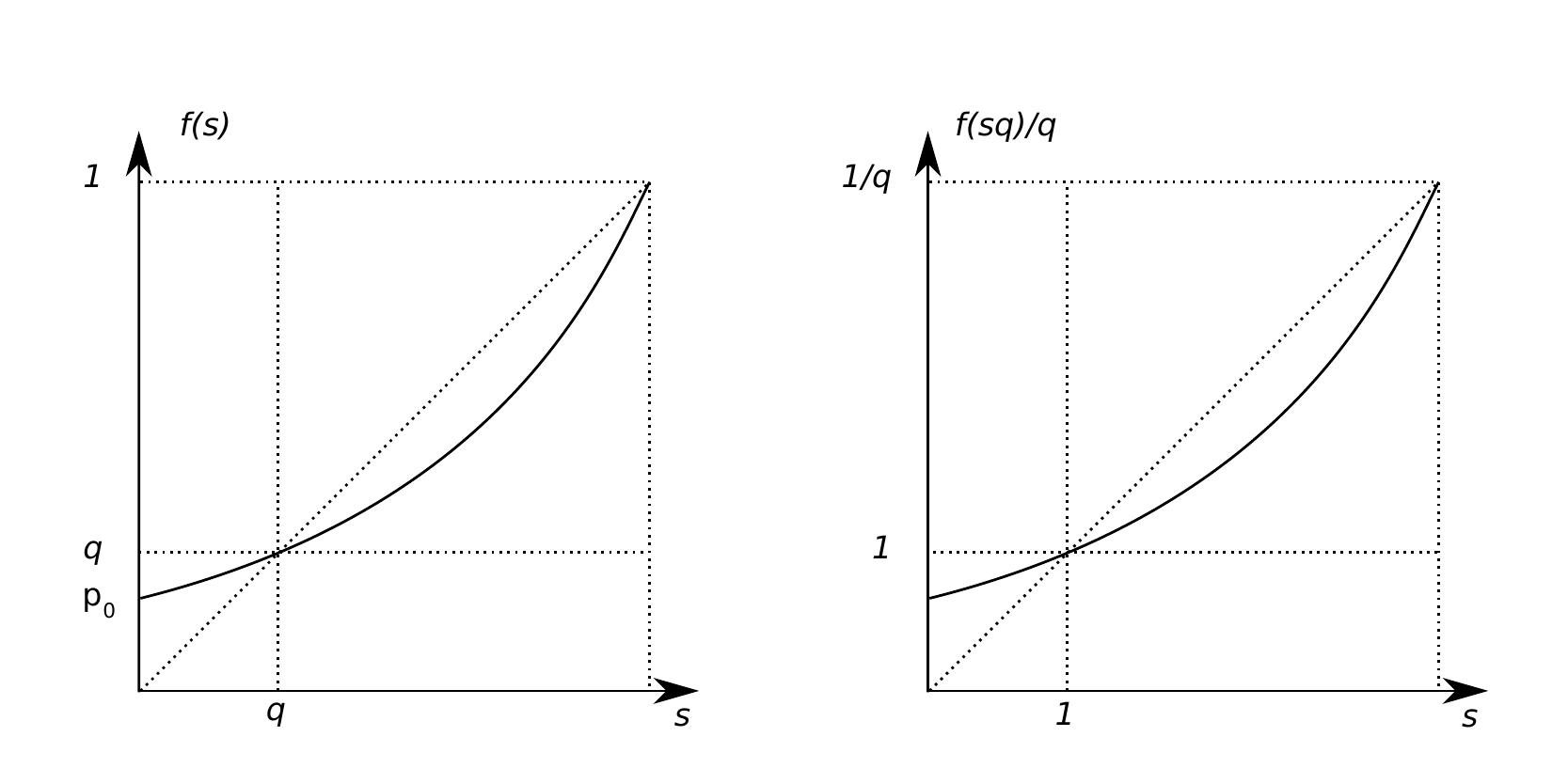}
\caption{Duality between the subcritical and supercritical cases. Left: a supercritical generating function \eqref{fs} with two positive roots $(q,1)$ for the equation \eqref{fx}. Right: the dual generating function  drawn on a different scale. }
\label{fig2}
\end{figure}

The corresponding dual BGW-process  is a subcritical branching process with offspring mean $\hat M=f'(q)<1$, see Figure \ref{fig2}.
The dual BGW-process is distributed as the original supercritical BGW-process conditioned on extinction:
\begin{eqnarray*}
 \mathbb P(Z^{(n)}=j|Z^{(0)}=i,Z^{(\infty)}=0) &=&{\mathbb P(Z^{(n)}=j,Z^{(\infty)}=0|Z^{(0)}=i)\over\mathbb P(Z^{(\infty)}=0|Z^{(0)}=i)}  \\
 &=&q^{-i}\mathbb P(Z^{(\infty)}=0|Z^{(n)}=j)P^{(n)}_{ij}  \\
 &=& q^{j-i}P^{(n)}_{ij}\\
 & =& \hat P^{(n)}_{ij}.
\end{eqnarray*}
%\hspace{1.5cm}
%\fbox{\parbox[c]{0.6\linewidth}{\small Supercritical multitype general BP conditioned on extinction is a subcritical BP. \\ Jagers and Lager\aa s (2008)
%}}
%
%\vspace{-2cm}\includegraphics[height=2cm]{jagers.jpg}
%\hspace{7cm}\vspace{-2cm}\includegraphics[height=1.8cm]{lageras.png}

The two parts of the curve on the left panel of Figure \ref{fig2} represent two transformations of the supercritical branching process. The lower-left  part of the curve, replicated on the right panel of Figure \ref{fig2} using a different scale, gives the generating function of the  dual process. The upper-right of the curve on the left panel corresponds to the  Harris-Sevastyanov transformation \eqref{HST}.
The function \eqref{HST} is the  generating function for the probability distribution
\[\tilde p_0=0, \ \tilde p_k=\sum_{i=k}^\infty p_i {i\choose k}q^{i-k}(1-q)^{k-1}\]
with the same mean $\tilde M=M$ as the original offspring distribution. It is easy to see that the $n$-th iteration of $\tilde f(s)$ is given by
\begin{align*}
\tilde f^{(n)}(s)&={f^{(n)}(s(1-q)+q)-q\over 1-q}.
\end{align*}

Looking into the future of the system of reproducing particles we can distinguish between two subtypes of particles:
\begin{quote}
\begin{itemize}
\item skeleton particles with infinite line of descent (building the skeleton of the genealogical tree),
\item doomed particles having finite line of descent.
\end{itemize}
\end{quote}
These two subtypes form a decomposable two-type BGW-process $\{(S^{(n)},D^{(n)})\}_{n=0}^\infty$ with $S^{(n)}+D^{(n)}=Z^{(n)}$. The joint reproduction law for the skeleton particles has the following generating function
 \begin{align*}
F(s,t)=\mathbb E(s^{S^{(1)}}t^{D^{(1)}})={f(s(1-q)+tq)-f(tq)\over 1-q}.
\end{align*}
A check on the branching property for the decomposed process is given by
\begin{align*}
\mathbb E(s^{S^{(n)}}t^{D^{(n)}}|Z^{(\infty)}>0)&={\mathbb E(s^{Z^{(n)}_1}t^{Z^{(n)}_2}1_{\{Z^{(\infty)}>0\}})\over \mathbb P(Z^{(\infty)}>0)}\\
&={\mathbb E(s^{S^{(n)}}t^{D^{(n)}})-\mathbb E(t^{D^{(n)}}1_{\{Z^{(\infty)}=0\}})\over1-q}\\
&={f^{(n)}(s(1-q)+tq)-f^{(n)}(tq)\over 1-q}=F^{(n)}(s,t).
\end{align*}
The original offspring distribution can be recovered  as a mixture of the joint reproduction laws of the two subtypes
 \begin{align*}
 f^{(n)}(s)&=(1-q)F^{(n)}(s,s)+q\hat f^{(n)}(s).
\end{align*}
Observe also that the total number of offspring for a skeleton particle has a distribution given by
\[\bar f(s)=F(s,s)={f(s)-f(sq)\over 1-q}=\sum_{k=0}^\infty \bar p_ks^k,\quad \bar p_k=(1-q)^{-1}(1-q^{k})p_k,\]
with mean $\bar M={M-\hat Mq\over1-q}$. It follows, 
\[\bar M=M+(M-\hat M){q\over1-q}\]
and we can summarize the relationship among different offspring means as 
$$\hat M<1<\tilde M=M<\bar M.$$

\section{ Linear-fractional single-type BGW-process}\label{sLF}
An important example of BGW-processes is the linear-fractional branching process. Its reproduction law has a  linear-fractional generating fuction
\begin{align}\label{lff}
f(s)=h_0+{h_1s\over1+m-ms}
\end{align}
fully characterized by two parameters: the probability $h_0=p_0$ of having no offspring, and the mean $m$ of the geometric number of offspring beyond the first one. Here $h_1=1-h_0$ stands for the probability of having at least one offspring. 
Notice that with 
$h_0={1\over 1+m}$, the generating function \eqref{lff} describes a Geometric $({1\over 1+m})$ distribution with mean $m$. 
If 
$h_0=0$ the generating function \eqref{lff} gives a Shifted Geometric $({1\over 1+m})$ distribution with mean $m+1$. 
If $m=0$, we arrive at a Bernoulli $(h_1)$ distribution.

Since the iterations of the linear-fractional function are again linear-fractional, many key characteristics of the linear-fractional BGW-processes can be computed explicitly in terms of the parameters $(h_0,m)$.
For example, we have $M=h_1(1+m)$, and if $M>1$, we get 
$$q=h_0(1+m^{-1})={1+m-M\over m}.$$

The dual reproduction law for \eqref{lff} is again linear-fractional
\begin{align*}
\hat f(s)&=\hat h_0+{\hat h_1s\over1+\hat m-\hat ms}, \qquad \hat h_0={m\over m+1}, \qquad \hat m=h_0/h_1,
\end{align*}
with $\hat M=1/M$.
The Harris-Sevastyanov transformation in the linear-fractional case corresponds to a shifted geometric distribution
\begin{align*}
\tilde f(s)&={s\over1+\tilde m-\tilde ms},\hspace{7mm} \tilde m=m(1-q)=M-1.
\end{align*}

Interestingly, the joint reproduction law of skeleton particles
\begin{align*}
    F(s,t)    &={h_1\over 1-q}\left({s(1-q)+tq\over 1+m-m(s(1-q)+tq)}-{tq\over 1+m-mtq}\right)\\
&={s\over1+m-m(s(1-q)+tq)}\cdot{1\over1+\hat m-\hat mt}
\\
&=s\cdot{1\over1+m-\tilde ms-(m-\tilde m)t}\cdot{1\over1+\hat m-\hat mt}
\end{align*}
has three independent components:
\begin{quote}
\begin{itemize}
\item one particle of type 1 (the infinite lineage),
\item a Geometric $({1\over 1+m})$ number of offspring each choosing independently between the skeleton and doomed subtypes with probabilities $\tilde m/m$ and $(m-\tilde m)/m$,
\item a Geometric $({1\over 1+\hat m})$ number of doomed offspring.
\end{itemize}
\end{quote}
Observe that even though both marginal distributions $\tilde f(s)$ and $\hat f(s)$
are linear-fractional, the decomposable BGW-process $(S^{(n)},D^{(n)})$ is not a two-type linear-fractional BGW-process. The distribution of the total number of offspring for the skeleton particles in not linear-fractional
\begin{align*}
   \bar f(s)  ={s\over1+m-ms}\cdot{1\over1+\hat m-\hat ms}
\end{align*}
and has mean 
$$\bar M=1+m+\hat m=M+(M+1)\hat m.$$

\section{BGW-processes  with countably many types}\label{sC}
A  BGW-process with countably many types
\[\mathbf{Z}^{(n)}=(Z_1^{(n)},Z_2^{(n)},\ldots),\ n=0,1,2,\ldots\]
describes demographic changes in a population of particles with different reproduction laws depending on the type of a particle $i\in\{1,2,\ldots\}$. Here $Z_i^{(n)}$ is the number of particles of type $i$ existing at generation $n$.
In the multi-type setting we use the following vector notation:
 \begin{align*}
\mathbf x&=(x_1,x_2\ldots), \ \mathbf 1=(1,1\ldots), \ \mathbf{e}_i=(1_{\{i=1\}},1_{\{i=2\}},\ldots),  \\
\mathbf x\mathbf y&=(x_1y_1,x_2y_2,\ldots), \hspace{3mm}\mathbf x^{-1}=(x_1^{-1},x_2^{-1},\ldots), \hspace{3mm} \mathbf x^{\mathbf y}=x_1^{y_1}x_2^{y_2}\cdots,
\end{align*}
we write $\mathbf x^{\rm t}$, if we need a column version of a vector $\mathbf x$.

A particle of type $i$ may produce random numbers of particles of different types so that the corresponding joint reproduction laws are given by the multivariate generating functions
\be
f_i(\mathbf{s})=\mathbb{E}( \mathbf{s}^{\mathbf{Z}^{(1)}}|\mathbf{Z}^{(0)}=\mathbf{e}_i).
\label{fis}
\ee
The offspring means
$$M_{ij}=\mathbb{E}(Z^{(1)}_j|\mathbf{Z}^{(0)}=\mathbf{e}_i)$$
 are convenient to summarize in a matrix form  $\mathbf{M}=\left(M_{ij}\right)_{i,j=1}^\infty$.
 For the $n$-th generation the vector of generating functions $\mathbf{f}^{(n)}(\mathbf{s})$ with components
$$f_i^{(n)}(\mathbf{s})=\mathbb{E}( \mathbf{s}^{\mathbf{Z}^{(n)}}|\mathbf{Z}^{(0)}=\mathbf{e}_i)$$
are obtained as iterations of $\mathbf{f}(\mathbf{s})$ with components \eqref{fis}, and the matrix of means is given by $\mathbf{M}^n$.
The vector of extinction probabilities $\mathbf{q}=(q_1,q_2,\ldots)$ has its $i$-th component $q_i$ defined as the probability of extinction given that the BGW-process starts from a particle of type $i$. The vector $\mathbf{q}$ is found as the minimal solution with non-negative components of equation $\mathbf{f}(\mathbf{x})=\mathbf{x}$, which is a multidimensional version of \eqref{fx}.

From now on we restrict our attention  to the positive recurrent (with respect to the type space) case when there exists a Perron-Frobenius eigenvalue $\rho$ for $\mathbf{M}$  with positive eigenvectors $\mathbf{u}$ and $\mathbf{v}$ such that
$$\mathbf{v}\mathbf{M}=\rho\mathbf{v}, \ \ \mathbf{M}\mathbf{u}^{\rm t}=\rho\mathbf{u}^{\rm t}, \ \ \mathbf v\mathbf u^{\rm t}=\mathbf v\mathbf 1^{\rm t}=1,  \ \ \rho^{-n}\mathbf{M}^n\to \mathbf{u}^{\rm t}\mathbf{v}, \ n\to\infty.$$
In the supercritical case, $\rho>1$, all $q_i<1$ and we can speak about the decomposition
 of a supercritical BGW-process with countably many types:
 $(\mathbf{S}^{(n)},\mathbf{D}^{(n)})$. Now each type is decomposed in two subtypes: either with infinite or finite line of descent. The decomposed supercritical BGW-process is again a BGW-process with countably many types whose reproduction law is given by the expressions
 \begin{align*}
F_{i}(\mathbf s,\mathbf t)={f_i(\mathbf s(\mathbf 1-\mathbf q)+\mathbf t\mathbf q)-f_i(\mathbf t\mathbf q)\over 1-q_i}, \ \ \ \hat f_{i}(\mathbf t)&={f_i(\mathbf t\mathbf q)\over q_i}.
\end{align*}

Linear-fractional BGW-processes with countably many types were studied recently in \cite{S}. In this case the joint probability generating functions \eqref{fis} have a restricted linear-fractional form
\be\label{fislf}
f_i(\mathbf s)=h_{i0}+{\sum_{j=1}^\infty h_{ij}s_j\over 1+m-m\sum_{j=1}^\infty g_js_j}.
\ee
The defining parameters of this branching process form a  triplet $(\mathbf H,\mathbf g,m)$, where $\mathbf H=(h_{ij})_{i,j=1}^\infty$ is a sub-stochastic matrix,  
$\mathbf g=(g_1,g_2,\ldots)$ is a proper probability distribution, and $m$ is a positive constant. The free term in \eqref{fislf} is defined  as
\[h_{i0}=1-\sum_{j=1}^\infty h_{ij}.\]
As in the case of finitely many types \cite {JL}, the denominators in \eqref{fislf} are necessarily independent of the mother type to ensure that the iterations are also linear-fractional. This is a major restriction of the multitype linear-fractional BGW-process excluding for example decomposable branching processes.

It is shown in \cite{S} that in the linear-fractional case the Perron-Frobenius eigenvalue $\rho$, if exists, is the unique positive solution of the equation
\be
m\sum_{k=1}^\infty \rho^{-k}\mathbf{g}\mathbf H^k\mathbf{1}^{\rm t}=1.
\label{PF}
\ee
In the positive recurrent case, when  the next sum is finite
\be
\beta=m\sum_{k=1}^\infty k \rho^{-k}\mathbf{g}\mathbf H^k\mathbf{1}^{\rm t},
\label{beta}
\ee
the Perron-Frobenius eigenvectors $(\mathbf v,\mathbf u)$ can be normalized in such a way that $\mathbf v\mathbf u^{\rm t}=\mathbf v\mathbf 1^{\rm t}=1$. They are computed as
\begin{align}
\mathbf u^{\rm t}&=(1+m)\beta^{-1}\sum_{k=1}^\infty \rho^{-k}\mathbf H^k\mathbf{1}^{\rm t},\label{ut}\\
\mathbf v&={m\over 1+m}\sum_{k=0}^\infty \rho^{-k}\mathbf{g}\mathbf H^k.\label{vt}
\end{align}
In the supercritical positive recurrent case with $\rho>1$ and  $\beta<\infty$ the extinction probabilities are given by
\begin{align}
\mathbf{q}&=\mathbf{1}-(\rho-1)(1+m)^{-1}\beta\mathbf{u}.
\label{gq}
\end{align}
Observe that $\mathbf g\mathbf u^{\rm t}={1+m\over m\beta}$ and
\begin{align}
\mathbf{g}\mathbf{q}^{\rm t}={1+m-\rho\over m}.\label{geq}
\end{align}
The total offspring number for a type $i$ particle has mean
\be
M_i=(1-h_{i0})(1+m).\label{Mi}
\ee

\section{ Main results}\label{sMR}

 In this section we summarize explicit formulae that we were able to obtain for the decomposition of the supercritical linear-fractional BGW-processes with countably many types. The derivation of these results is given in the next section.
 
 Consider the positive recurrent supercritical case with $\rho>1$ and  $\beta<\infty$.  We demonstrate that the dual reproduction laws are again linear-fractional
\be
\hat f_i(\mathbf s)=\hat h_{i0}+{\sum_{j=1}^\infty \hat h_{ij}s_j\over 1+\hat m-\hat m\sum_{j=1}^\infty \hat g_js_j},
\label{fi}
\ee
with
\begin{align}
\hat h_{i0}&={h_{i0}\over q_i}, \hspace{19mm}\hat h_{ij}={h_{ij}q_j\over q_i\rho}, \label{fi1}\\
\hat m&={1+m-\rho\over \rho}, \hspace{7mm} \hat g_{j}={g_{j}q_jm\over 1+m-\rho}.\label{fi2}
\end{align}
%Notice that in view of \eqref{Q} and $\rho>1$ we get $\hat m<m$.
It turns out that the following remarkably simple formulae hold for the key characteristics of  the dual branching process
\begin{align}
\hat \rho&=\rho^{-1},  \label{hr}\\
 \hat \beta&={\mu-1\over \rho-1},\quad \mu=m\sum_{n=1}^\infty \mathbf{g}\mathbf H^n\mathbf{1}^{\rm t}. \label{br}
\end{align}
For the Perron-Frobenius eigenvectors we obtain the following expressions
\begin{align}
\hat {\mathbf u}&=\beta \hat \beta^{-1}\mathbf u\mathbf q^{-1}=(\mathbf q^{-1}-\mathbf 1)(1+m)( \mu-1)^{-1},\label{u}\\
 \hat {\mathbf v}&=\frac{m}{1+m}\sum_{k=0}^\infty(\mathbf{g}\mathbf{H}^k)\mathbf q.\label{v}
\end{align}

We show that the Harris-Sevastyanov transformation results in multivariate shifted geometric distributions
\be
\tilde f_i(\mathbf s)={\sum_{j=1}^\infty \tilde  h_{ij}s_j\over 1+\tilde  m-\tilde  m\sum_{j=1}^\infty \tilde  g_js_j}, \label{tfi}
\ee
where 
\begin{align}
\tilde h_{ij}&={1-q_j\over 1-q_i}\left(h_{ij}+mg_j(q_i-h_{i0})\right),\label{mgh}\\
\tilde m&=\rho-1,\quad \tilde  g_{j}={m\over \rho-1}g_{j}(1-q_j). \label{mgh1}
\end{align}
Moreover, we demonstrate that
\be
\tilde  \rho=\rho,\qquad  \tilde  \beta={\rho\over \rho-1},\label{tb}
 \ee
and
\begin{align}
\tilde  {\mathbf u}=\mathbf 1,\hspace{5mm}\tilde{\mathbf{v}}
&=m\sum_{k=0}^\infty\rho^{-1-k}\Big(\mathbf{g}\{\mathbf{H}+m\rho^{-1}\mathbf{H}\mathbf{q}^{\rm t}\mathbf{g}\}^k\Big)(\mathbf{1}-\mathbf{q}). \label{uv}
\end{align}

\begin{theorem}\label{Thm}
 Consider a linear-fractional BGW-process characterized by a triplet $(\mathbf H,\mathbf g,m)$. Assume it is supercritical and positively recurrent over the state space, that is $\rho>1$ and $\beta<\infty$. 
 Its dual BGW-process and its skeleton are also linear-fractional BGW-processes with the transformed parameter triplets $(\hat{\mathbf H},\hat{\mathbf g},\hat m)$ and $(\tilde{\mathbf H},\tilde{\mathbf g},\tilde m)$ with components given by \eqref{fi1}, \eqref{fi2}, \eqref{mgh}, \eqref{mgh1}. 
 
 The joint offspring generating function for a skeleton particle of type $i$ has the form
\begin{align}
 F_{i}(\mathbf s,\mathbf t)&=\sum_{j=1}^\infty { \tilde h_{ij} s_j\over 1+m-\tilde m\tilde{\mathbf g}\mathbf s^{\rm t}-(m-\tilde m)\hat{\mathbf g}\mathbf t^{\rm t}}\Big(h_{ij0}+{\sum_{k=1}^\infty  h_{ijk}t_k\over 1+\hat m- \hat m\hat{\mathbf g}\mathbf t^{\rm t}}\Big),
\label{f1i}
\end{align}
where
\[h_{ij0}={h_{ij}\over h_{ij}+mg_j(q_i-h_{i0})},\quad h_{ijk}={mg_jq_i  \hat h_{ik}\over h_{ij}+mg_j(q_i-h_{i0})}.\]
\end{theorem}
Similarly to the single-type case, we can distinguish in \eqref{f1i} three components but now with dependence:
\begin{itemize}
\item a ``reborn" skeleton particle of type $i$ may change its type to $j$ with probability $\tilde h_{ij}$,
\item independent of $i$ and $j$ a multivariate geometric number of offspring of both subtypes,
\item a linear-fractional number of doomed offspring with the fate of the first offspring being dependent on $(i,j)$ .
\end{itemize}
The total number of offspring of a skeleton particle of type $i$ has generating function $ \bar f_i(s)=F_{i}(s\mathbf 1,s\mathbf 1)$ of the next form
\begin{align*}
 \bar f_i(s)%&={s\over 1+m-ms}\sum_{j=1}^\infty{1-q_j\over 1-q_i}\left(h_{ij}+{mg_jq_i(1-\hat h_{i0})s\over 1+\hat m- \hat ms}\right)\\
&={s\over 1+m-ms}\Big(1-\alpha_i+{\alpha_is\over 1+\hat m- \hat ms}\Big),
\end{align*}
where $\alpha_i={\rho-1\over 1-q_i}(q_i-h_{i0})$ must belong to the interval $(0,1)$. The corresponding mean offspring number is larger than that given by \eqref{Mi}:
\begin{align*}
\bar M_i&=1+m+\alpha_i(1+\hat m)=M_i+(1+m)\Big(h_{i0}+{(\rho-1)(q_i-h_{i0})\over \rho(1-q_i)}\Big).
\end{align*}

%The corresponding decomposition of the associated linear-fractional single-type CMJ-processes should be viewed in terms of stem individuals and leaf individuals. The stem individuals have infinite life lengths. The life length under the dual law
%\begin{align}
%\hat {\mathbb P}({L}=n)&=\frac{m}{(1+m-\rho)\rho^{n-1}}\mathbb P(L=n)
%\label{Ln}
%\end{align}
%is obtained as a change of measure transformation with more weight allocated to smaller values. It follows that 
%\be
%\hat {\mathbb E}({L})=1+\frac{\rho(1+m-\beta(\rho-1))}{(1+m-\rho)^2}
%\label{EL}
%\ee
%and
%\begin{align}
%\hat \mu={1+m-\beta(\rho-1)\over 1+m-\rho}.
%\label{hmu}
%\end{align}
%Observe that since $\hat \mu<1$, we have $\beta>{\rho\over \rho-1}$.
%Interestingly, the mean age at childbearing $ \tilde  \beta$
%is finite even if $\mathbb {\tilde P}(L=\infty)=1$. The infinite life length implies $\tilde \mu=\infty$.
%
%Clearly, $\hat Q=1$ and $\tilde Q=0$ and the initial individual of the decomposed linear-fractional CMJ-process is a stem individual with probability $1-Q$ and a leaf individual with  probability $Q$ given by \eqref{Q}.

\section{ Proof of Theorem \ref{Thm}}\label{sP}

In this section we derive the formulae stated in Section \ref{sMR}.

\vspace{5mm}{\bf Proof of \eqref{fi}.}
From
 \[\hat f_i(\mathbf s)={f_i(\mathbf s\mathbf q)\over q_i}=\frac{h_{i0}}{q_i}+\frac{\sum_{j=1}^{\infty}h_{ij}s_j q_j/q_i}{1+m-m \sum_{j=1}^{\infty} g_j s_j q_j}\]
it is straightforward to obtain \eqref{fi} with \eqref{fi1} and \eqref{fi2}. We have to verify that $\hat{\mathbf g}\mathbf{1}^{\rm t}=1$ and
\[\hat h_{i0}=1-\sum_{j=1}^\infty\hat h_{ij}.\]
The first requirement follows from 
 \eqref{geq}. The second is obtained from
\begin{align}
\mathbf{H}\mathbf{q}^{\rm t}
&=\rho(\mathbf{H}\mathbf{1}^{\rm t}-\mathbf{1}^{\rm t}+\mathbf{q}^{\rm t}) \label{Hq}
\end{align}
which is proved next. We have (relation (6) in \cite{S})
\[\mathbf H=\mathbf M-{m\over 1+m}\mathbf M\mathbf{1}^{\rm t}\mathbf g\]
and therefore $\mathbf M\mathbf{1}^{\rm t}=(1+m)\mathbf H\mathbf{1}^{\rm t}$, which is  \eqref{Mi}.
Using the last two equalities and \eqref{gq} we find first
\begin{align}
\mathbf{H}\mathbf{q}^{\rm t}
&={1\over 1+m}\mathbf{M}\mathbf{1}^{\rm t}-\rho(\mathbf{1}^{\rm t}-\mathbf{q}^{\rm t}) +{\rho-1\over 1+m}\mathbf{M}\mathbf{1}^{\rm t}\nonumber 
\end{align}
and then obtain \eqref{Hq}.

\vspace{5mm}{\bf Proof of  \eqref{hr}.} In view of equation \eqref{PF} determining the Perron-Frobenius eigenvalue for a linear-fractional BGW-process, to show \eqref{hr} it is enough to verify that
\begin{align*}
\hat{m} \sum_{n=1}^\infty \rho^{n} \hat{\mathbf{g}} \hat{\mathbf{H}}^n \mathbf{1}^{\rm t}=1.
\end{align*}
Observe that according to \eqref{fi1}
\be
\hat{\mathbf{g}}\hat{\mathbf{H}}^k=\frac{m}{(1+m-\rho)\rho^k}(\mathbf{g}\mathbf{H}^k) \mathbf{q}.
\label{gH}
\ee
It follows,
\begin{align}
\hat{\mathbf{g}} \hat{\mathbf{H}}^{n} \mathbf{1}^{\rm t}=\frac{m}{(1+m-\rho)\rho^{n}} \mathbf{g}\mathbf{H}^{n}\mathbf{q}^{\rm t},
\label{gH1}
\end{align}
so that we have to check that
\begin{align}
m\sum_{n=1}^\infty \mathbf{g}\mathbf{H}^{n}\mathbf{q}^{\rm t}=\rho.\label{mur}
\end{align}
Turning to \eqref{Hq} we find
\begin{align}
\mathbf{H}^n\mathbf{q}^{\rm t}
&=\rho(\mathbf{H}^n\mathbf{1}^{\rm t}-\mathbf{H}^{n-1}\mathbf{1}^{\rm t}+\mathbf{H}^{n-1}\mathbf{q}^{\rm t}) \label{Hn}
\end{align}
yielding
\begin{align}
    (\rho-1)\sum_{n=1}^\infty \mathbf{H}^{n}\mathbf{q}^{\rm t} =\rho(\mathbf{1}^{\rm t}-\mathbf{q}^{\rm t}).\label{qur}
\end{align}
This and \eqref{geq} entail \eqref{mur}.

\vspace{5mm}{\bf Proof of  \eqref{br}.} Starting from a counterpart of \eqref{beta}  we find using \eqref{gH1}
\begin{align*}
\hat{\beta}&=\hat{m}\sum_{n=1}^\infty n
\hat{\rho}^{-n}\hat{\mathbf{g}}\hat{\mathbf
H}^n\mathbf{1}^{\rm t}=\frac{m}{\rho}\sum_{n=1}^\infty n \mathbf{g} \mathbf{H}^{n} \mathbf{q}^{\rm t}.
\end{align*}
Rewrite \eqref{Hn} as 
\begin{align*}
n\mathbf{H}^n\mathbf{q}^{\rm t}
&=\rho\Big(n\mathbf{H}^n\mathbf{1}^{\rm t}-(n-1)\mathbf{H}^{n-1}\mathbf{1}^{\rm t}+(n-1)\mathbf{H}^{n-1}\mathbf{q}^{\rm t}\Big)+\rho \mathbf{H}^{n-1}(\mathbf{q}^{\rm t}-\mathbf{1}^{\rm t})%\label{Hn}
\end{align*}
to obtain
$$\sum_{n=1}^\infty n\mathbf{H}^{n}\mathbf{q}^{\rm t}=\frac{\rho}{\rho-1}\sum_{n=0}^\infty \mathbf{H}^{n}(\mathbf{1}^{\rm t}-\mathbf{q}^{\rm t}).$$
Thus 
$$\hat{\beta}=\frac{m}{\rho-1}\sum_{n=0}^\infty \mathbf{g}\mathbf{H}^{n}(\mathbf{1}^{\rm t}-\mathbf{q}^{\rm t})=\frac{\mu-1}{\rho-1}.$$

\vspace{5mm}{\bf Proof of  \eqref{u} and \eqref{v}.} From \eqref{fi1} we derive $\hat{\mathbf{H}}^{n}\mathbf{1}^{\rm t}=(\mathbf{H}^{n}\mathbf{q}^{\rm t})\mathbf{q}^{-1}$.
This and a counterpart of \eqref{ut}
$$ \hat{\mathbf{u}}=(1+\hat{m}) \hat{\beta}^{-1}\sum_{n=1}^\infty \hat{\rho}^{-n}\hat{\mathbf{H}}^{n}\mathbf{1}^{\rm t}$$
in view of \eqref{qur} brings \eqref{u}
$$ \hat{\mathbf{u}}=
\frac{(1+m)(\rho-1)}{\rho(\mu-1) } \Big(\sum_{n=1}^\infty \mathbf{H}^{n}\mathbf{q}^{\rm t}\Big)\mathbf{q}^{-1}=(\mathbf q^{-1}-\mathbf 1)(1+m)( \mu-1)^{-1}.$$
On the other hand, a counterpart of \eqref{vt} together with \eqref{gH} yields
$$\hat{\mathbf{v}}=\frac{\hat{m}}{1+\hat{m}}
\sum_{k=0}^\infty \hat{\rho}^{-k}
\hat{\mathbf{g}}\hat{\mathbf{H}}^k=\frac{m}{1+m} \sum_{k=0}^\infty \mathbf{g}\mathbf{H}^k \mathbf{q}.$$ 

\vspace{5mm}{\bf Proof of  \eqref{f1i}.} We have
\begin{align*}
&f_i(\mathbf s(\mathbf 1-\mathbf q)+\mathbf t\mathbf q)-f_i(\mathbf t\mathbf q)\\
&={\sum_{j=1}^\infty h_{ij}s_j(1-q_j)+\sum_{j=1}^\infty h_{ij}t_jq_j\over 1+m-m\sum_{k=1}^\infty g_ks_k(1-q_k)-m\sum_{k=1}^\infty g_kt_kq_k}-{\sum_{j=1}^\infty h_{ij}t_jq_j\over 1+m- m\sum_{k=1}^\infty g_kt_kq_k}\\
&={\sum_{j=1}^\infty h_{ij}s_j(1-q_j)
\over 1+m-\tilde m\sum_{k=1}^\infty \tilde g_ks_k-(m-\tilde m)\sum_{k=1}^\infty \hat g_kt_k}\\
&+{\big(\sum_{j=1}^\infty h_{ij}t_jq_j\big)\big(m\sum_{k=1}^\infty g_ks_k(1-q_k)\big)\over\big(1+m-\tilde m\sum_{k=1}^\infty \tilde g_ks_k-(m-\tilde m)\sum_{k=1}^\infty \hat g_kt_k\big)\rho\big(1+\hat m- \hat m\sum_{k=1}^\infty \hat g_kt_k\big)}.
\end{align*}
Replacing the last numerator by $m\sum_{j=1}^\infty g_js_j(1-q_j)\sum_{k=1}^\infty h_{ik}t_kq_k$ and dividing the whole expression by $1-q_i$ we get
\begin{align*}
 F_{i}(\mathbf s,\mathbf t)&=\sum_{j=1}^\infty {s_j(1-q_j)(1-q_i)^{-1}
\over1+m-\tilde m\tilde{\mathbf g}\mathbf s^{\rm t}-(m-\tilde m)\hat{\mathbf g}\mathbf t^{\rm t}}\left(h_{ij}+{mg_jq_i\sum_{k=1}^\infty  \hat h_{ik}t_k\over 1+\hat m- \hat m\hat{\mathbf g}\mathbf t^{\rm t}}\right)
\end{align*}
and the relation  \eqref{f1i} follows.

\vspace{5mm}{\bf Proof of \eqref{tfi}, \eqref{tb}.} Putting $\mathbf t=\mathbf 1$ in \eqref{f1i} we arrive at \eqref{tfi} . Notice that according to definition \eqref{mgh} and relations  \eqref{geq},  \eqref{Hq} we have
\[\tilde{\mathbf g}\mathbf 1^{\rm t}=1,\qquad \tilde{\mathbf H}\mathbf 1^{\rm t}=\mathbf 1^{\rm t}.\]
Since $\tilde{\rho}$ is the unique positive solution of
$$\tilde{m}\sum_{n=1}^\infty \tilde{\rho}^{-n} \tilde{\mathbf{g}}\tilde{\mathbf{H}}^n \mathbf{1}^{\rm t}=1$$
and $\tilde{\mathbf{g}}\tilde{\mathbf{H}}^{n}\mathbf{1}^{\rm t}=1$ we derive $$(\rho-1)\sum_{n=1}^\infty \tilde{\rho}^{-n}=1.$$ 
Thus $\tilde{\rho}=\rho$ and
$$\tilde{\beta}=\tilde{m}\sum_{n=1}^\infty n \tilde{\rho}^{-n}\tilde{\mathbf{g}}\tilde{\mathbf{H}}^n\mathbf{1}^{\rm t}=(\rho-1)\sum_{n=1}^\infty n
\rho^{-n}=\frac{\rho}{\rho-1}.$$

\section*{Acknowledgments}

SS was supported by the Swedish Research Council grant 621-2010-5623. AS was supported by the Scientific Committee of  Kazakhstan's Ministry of Education and Science, grant 0732/GF 2012-14.

\end{document}